# STORY OF AN ARCHITECTURALLY SUGGESTIVE POLYHEDRON: FROM MEDIEVAL TRADE TO RENAISSANCE ART AND MODERN DESIGN

EUGENE A. KATZ

*Ben-Gurion University of the Negev, Sede Boqer Campus, 8499000, Israel*
*E-mail: keugene@bgu.ac.il*

**Abstract:** *We describe one of the balance weights dated to the Early Islamic Period from the Hecht Museum at the University of Haifa (Israel). Its polyhedral shape was attributed to a truncated elongated octagonal bipyramid. To our knowledge, the earliest Renaissance book containing the image of this polyhedron is 'La Pratica di Prospettiva' published in 1596 by Florentine architect and perspective artist Lorenzo Sirigatti. We described an outline of Sirigatti's life and the importance of his book for scientists and artists, Galileo Galilei among them. We depict examples of lamps and lanterns in European cities shaped as the truncated elongated octagonal bipyramid, as well as a drawing by Raphael with a lantern of a similar form. Finally, we discussed why the artist chose this particular polyhedron for his drawing.*

## INTRODUCTION

Archaeologists commonly found the copper-alloy polyhedral balance weights in the Middle East and even across the Scandinavian world, from Russia in the east to Ireland and England in the west. The archaeological literature has a broad consensus that these findings reflect trading contacts with the Islamic caliphate, and most weights were imported from the Middle East during the VIII – X centuries (Kershaw, 2019). These artifacts are often called the dirham weights since they were used to weigh silver coins, *dirhams*.

Recently, studying such archeologic artifacts in the context of the history of their geometric forms (Chernov, 2024), I visited the Hecht Museum at the University of Haifa and found several polyhedral dirham weights dated to the Early Islamic Period (Fig. 1). I was allowed to photograph them, but the museum workers said the weights originated from a private purchase, and it is not possible to know for sure where they were found.

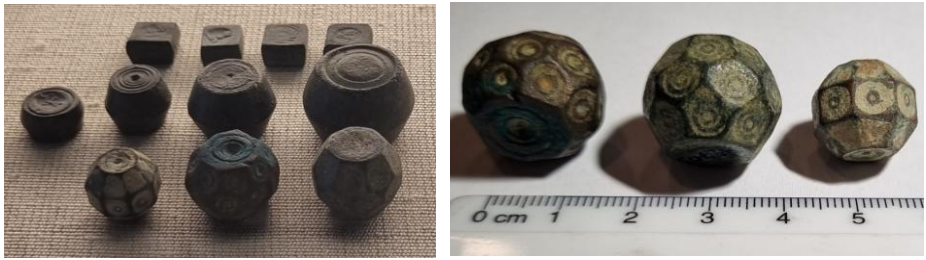

Fig. 1 The dirham weights from the Hecht Museum at the University of Haifa, Israel, dated to the Early Islamic Period (VIII – X centuries). Photo by E. A. Katz.

# 1 DIRHAM WEIGHT IN THE FORM OF TRUNCATED ELONGATED OCTAGONAL BIPYRAMID

Below, we will describe one of these polyhedra (Fig. 2) and investigate its appearance in Renaissance mathematical literature and art. This is a Truncated Elongated Octagonal Bipyramid (TEOB). It comprises a central row of eight rectangular (almost square) faces parallel to the vertical eight-fold axis of rotational symmetry, sixteen identical trapezium faces above and below the central raw (8+8), and two regular octagonal horizontal faces. In other words, it consists of the octagonal prism and the truncated octagonal pyramids above and below it. The TEOB can be described as "t8k8P8" in Conway notation (Conway, 2008) and generated with software https://levskaya.github.io/polyhedronisme/?recipe=t8k8P8.

The weight polyhedron can be approximately inscribed in a sphere with a diameter of ~1.5 cm. Tekin described a collection of balance weights at the Anatolian Civilization Museum in Ankara, Turkey (Tekin, 2016). Most of them date to the Early Islamic Period. Some weights may have TEOB form, though their images are unclear, and no description is given regarding their geometry.

For someone to produce weights with such sophisticated shapes, their geometry should be already known. In our previous study (Chernov, 2024), one of the arguments about the Middle Eastern origin of the dirham weights in the form of cuboctahedron was that the mathematicians of the Islamic Caliphate were aware of the cuboctahedron. Still, European mathematicians of that time were not (Field, 1997). In particular, all Platonic solids and some Archimedean solids were described by the Persian mathematician and astronomer Abū al-Wafā Būzhjānī (940–998) in the treatise "A Book on Those Geometric Constructions Which Are Necessary for a Craftsman" ( كتاب في ما يحتاج إليه الصانع من الأعمال الهندسية) written after 990 (Sarhangi, 2008). However, to our knowledge, the TEOB was described neither in this treatise nor in the famous Renaissance book "De Divina Proportione" (Pacioli, 1509) written by Luca Pacioli (1447–1517), as well as in the subsequent essential publications on the theory of perspective with extensive collections of polyhedra images (Barbaro, 1568; Jamnitzer, 1568).

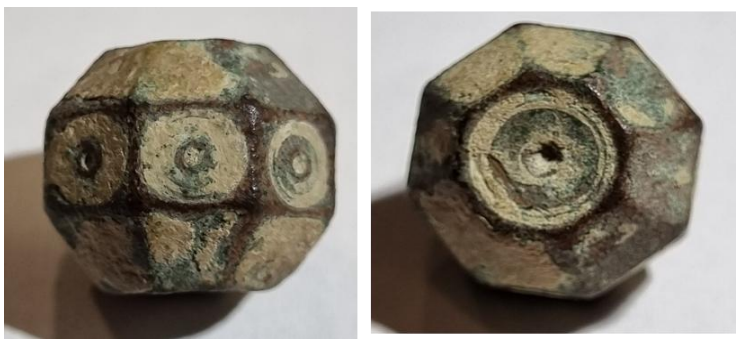

Fig 2. The dirham weight in the form of truncated elongated octagonal bipyramid from the Hecht Museum at the University of Haifa, Israel. Photo by E. A. Katz.

## 2 LORENZO SIRIGATTI AND HIS BOOK "LA PRATICA DI PROSPETTIVA"

To our best knowledge, the earliest Renaissance book containing the TEOB image is 'La Pratica di Prospettiva' by Florentine architect and perspective artist Lorenzo Sirigatti (1557- 1625), printed in 1596 in Venice (Sirigatti, 1596).

Little is known about the author's life (Pegazzano, 2018). He was born in Florence in a family with a long artistic tradition. His oldest brother, Ridolfo Sirigatti (1553–1608), was a sculpturer and painter, and his mother, Cassandra, was a granddaughter of the great Renaissance artist Domenico Ghirlandaio (1448 - 1494).

The original Sirigatti's manuscript, dated 1593, complete with his drawings, is preserved in the Laurentian Library in Florence (Ashburnham code 1029). The title of the treatise completely repeats that of the book by Venetian humanist and mathematician Daniele Barbaro (1514-1570) (Barbaro, 1568). This coincidence is probably not accidental and can be explained by Sirigatti's stay in Venice at a young age.

The book contains sixty-three engravings deriving from the original Sirigatti's drawings and depicting polyhedra (including all platonic solids, some Archimedean polyhedra, and the Campanus sphere), architectural elements, facades of palaces and churches, musical instruments. It is divided into two parts. In the first part, the engravings are accompanied by perspective exercises and the elementary rules of perspective. In particular, Sirigatti contributed to the study of theatrical perspective.

The second part is only a collection of the engravings. These images' clarity, the engraving's excellent quality, and a concise theory description made the treatise a practical manual used by generations of artists and architects.

The TEOB is shown only in the second part of the book (Fig. 3, left). It is shown inscribed into a cube, while the Campanus sphere is inscribed inside it. The TEOB and the cube are depicted in the "skeletal" manner ("basium vacuum") developed by Leonardo da Vinci for illustration of Luca Pacioli's "De Divina Proportione". There is no textual description of "our" polyhedron, but the detailed shadow pattern provides excellent instruction for artists or architects studying the perspective.

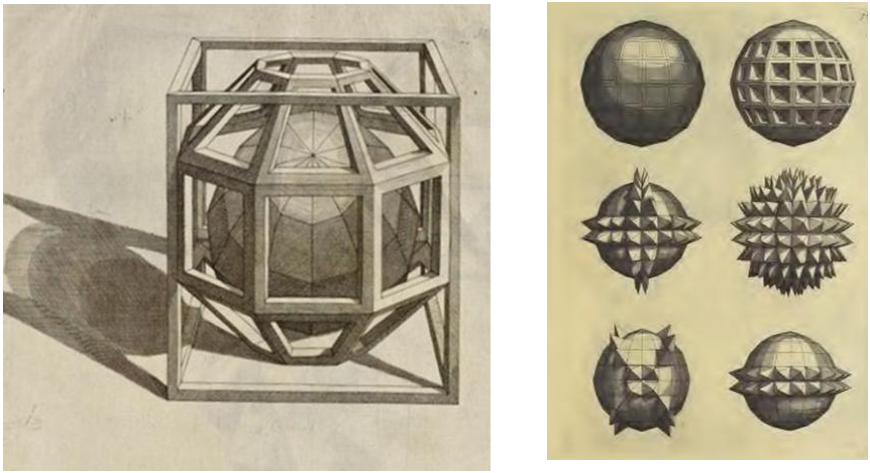

Fig. 3. Polyhedral sketches from Lorenzo Sirigatti's 'La pratica di prospettiva'. Left: the TEOB with the cube and the Campanus sphere. Right: Images of complex figures with detailed light/shadow patterns. Reproduced from https://archive.org/details/gri_33125009324340

Many were influenced by Sirigatti's book lessons on chiaroscuro, a technique for representing three-dimensional subjects on a plane using light/shadow patterns of complex geometric forms (Figure 4). The most famous was … Galileo Galilei (1564-1642).

Galileo studied Sirigatti's book, published in Venice, while Galileo was a professor of mathematics nearby at the University of Padova. Furthermore, Sirigatti was a member of the Academy of Drawing (Accademia delle Arti del Disegno), a school for painters, sculptors, and architects, where young Galileo had got his background in art, perspective, and particularly in the practice of chiaroscuro. Like any other student, he had to learn to show a pattern of light and shadow from the spikes on a ring diagram (Figure 4). Each spike must cast an appropriate shadow, not so unlike the patches Galileo would later observe using his telescope and interpret as the shadows of mountains protruding up from the surface of the Moon. Indeed, when Galileo and English astronomer Thomas Harriot (c. 1560 – 621) almost simultaneously pioneered the use of the telescope to study the moon's surface, it was Galileo's training in chiaroscuro that led him to see mountains and craters where Harriot only saw "strange spottiness" (Reeves, 1997).

In the years following the publication of 'La Pratica di Prospettiva' Sirigatti performed various architectural projects, most of which did not materialize. The only surviving Sirigatti's project is the facade of the Florentine palace of the Della Fonte family (today Sebregondi) in Florence (via Ghibellina, 81). The style of this facade shows faithful adherence to contemporary Florentine architecture (Fig. 4).

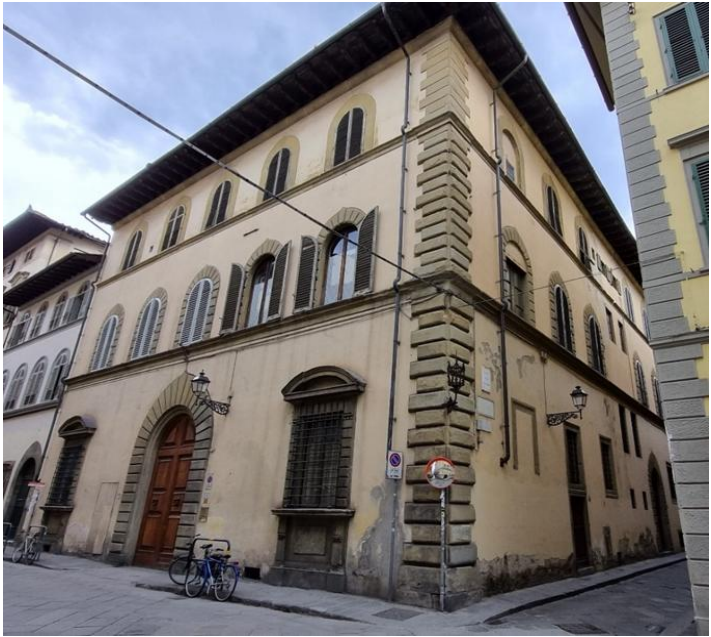

Fig. 4. Palazzo Salvetti Sebregondi, Florence, via Ghibellina, 81. Although having more ancient origins, the palace currently presents characteristics based on the project by Lorenzo Sirigatti , who unified and redesigned the pre-existing architecture commissioned by the Della Fonte family. Photo by E. A. Katz.

## 3 TEOB IN MODERN DESIGN

Lamp designers of modern times recognized the beauty of the TEOB. Lamps and lanterns with a TEOB-like cover (Fig. 5) are widely seen inside and outside of buildings/palaces in European cities. This tradition probably started in the Baroque. I took the photos in Fig. 5 (a) in Turin in Palazzo Carignano. This is a masterpiece of the Piedmontese Baroque designed in 1679 by the great architect and polymath Camillo Guarino Guarini (1624-1683). I am not sure that the original Guarini's design included this (or) similar lanterns. Further research is needed to check this suggestion. According to Edoardo Piccoli, an expert in Guarini's architecture, while the technique and decoration point to the late 19th century (or early 20th), we probably owe this lantern design to an architect emulating Guarini's elaborate style and geometry (Piccoli, 2024).

In a couple of days after getting the above-cited opinion, I received another letter from Prof. Piccoli about his exciting finding: "… while browsing my architectural course images, my eye fell on a drawing by Raphael, circa 1520, with the same lantern shape!". This Raphael's drawing is shown in Figure 6.

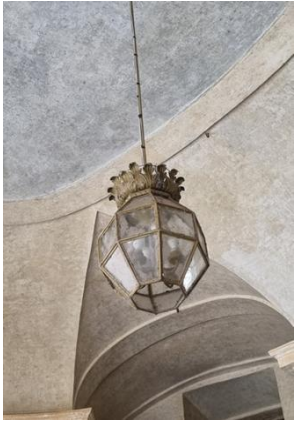 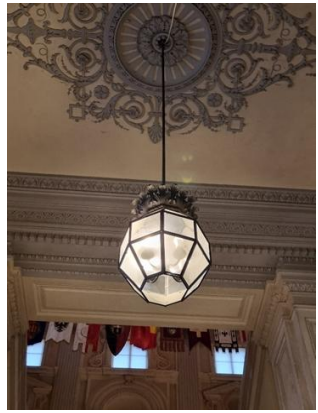
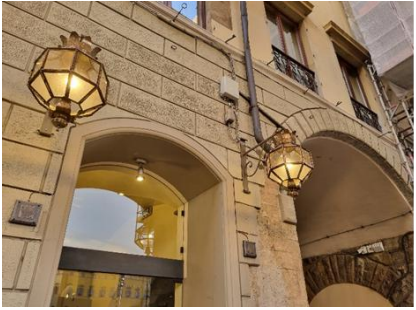 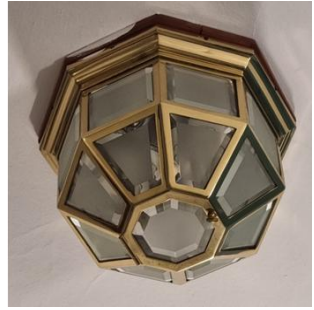

Fig. 5. Lamps with a cover shaped in the truncated elongated octagonal bipyramid sphere (a-c) or part of it (d). Photos made by E. A. Katz in Turino (a, b), Florence (c), and Bolzano (d).

## 4  FINDING THE RAPHAEL DRAWING RAISES NEW QUESTIONS

The by Raphael (1483-1520) does not correspond to any specific building and is classified as a perspective or architectural "fantasy". Attribution to Raphael is not controversial, even though the drawing is not signed. It is dated to the beginning of the 16th century (Wikipedia, 2024; Wikimedia, 2024). Raphael could see such a lamp and decide to use it as a light source in the center of the composition, whose real subject is an interplay of light and shade in the context of a linear perspective drawing.

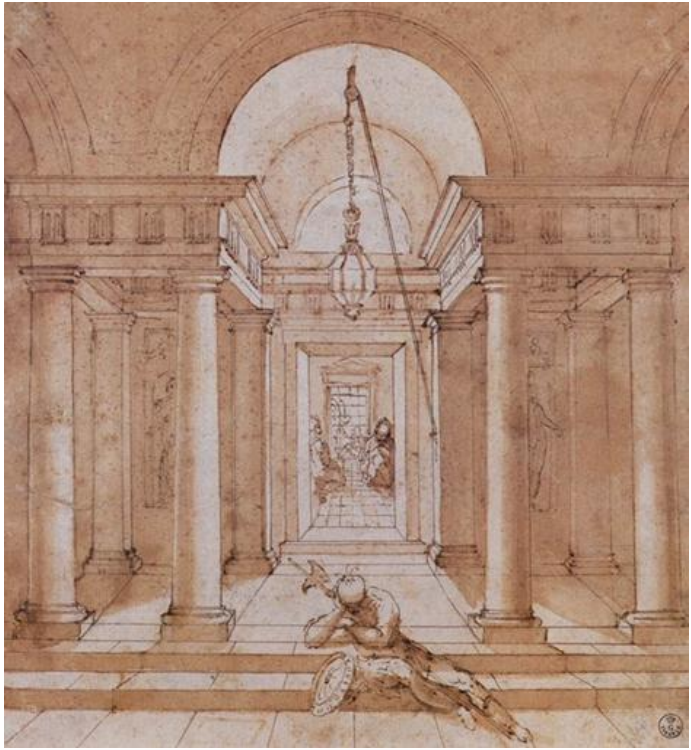

Fig. 6. Raphael, Soldier before the cell of St. Peter, 1515-1520, pen and ink and wash on white paper. Cabinet of Drawings and Prints of the Uffizi Gallery, Florence. Reproduced from Wikimedia, public domain (Wikimedia, 2024).

This finding is significant in our study since Raphael's drawing is dated much earlier than 1596 when Lorenzo Sirigatti printed his book. Though it does not challenge our conclusion that the first book depicting the truncated elongated octagonal bipyramid is Sirigatti's "La Pratica di Prospettiva", it raises several new questions.

What is known about the context of the Raphael drawing? Where could Raphael see such a polyhedron? Why did he choose this particular polyhedron? I could not find direct answers to any of these questions in art history literature. However, some specific help can be found in the paper "Archimedes Salutes Bramante in a Draft for the School of Athens" by Caroline Karpinski (Karpinski, 2010). It is focused on the analyses of another masterpiece of almost the same time, the chiaroscuro woodcut "Archimedes" (Fig. 7). It is attributed to 'Ugo da Carpi, after Raphael'. It is clear from the article title that the author believes that Ugo da Carpi reproduced in the engraving Raphael's drawing of Archimedes, as he was initially supposed to be depicted in Raphael's monumental fresco "The School of Athens" (1509-1511).

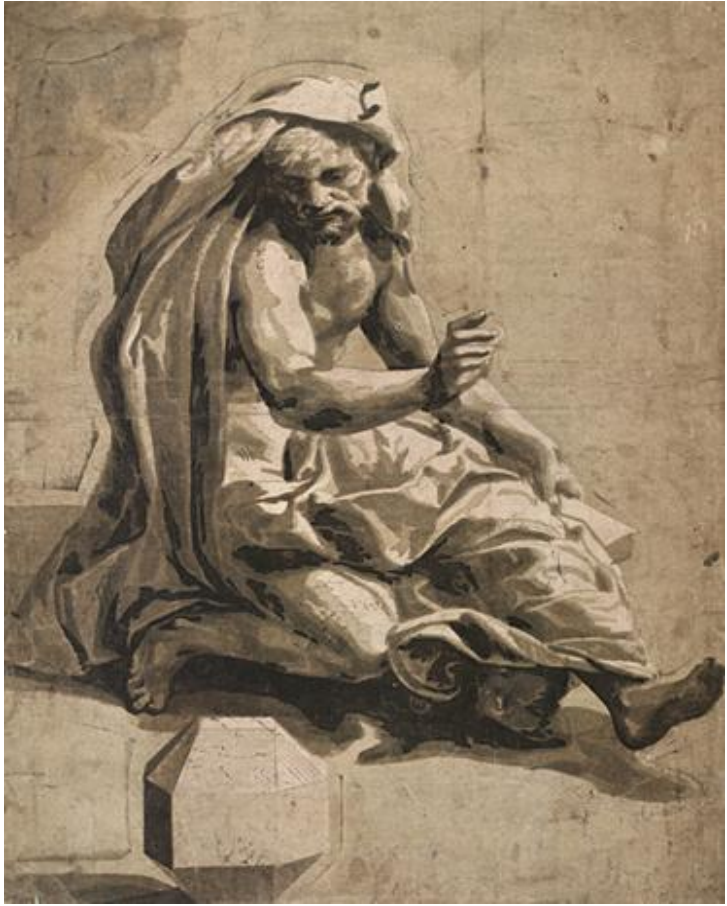

Fig. 7. Ugo da Carpi (c. 1450–1480 – c. 1523–1532), after Raphael. Archimedes. 1518-20, chiaroscuro woodcut. Albertina, Vienna. Reproduced from: https://timelessmoon.getarchive.net/amp/media/ugo-da-carpi-archimedes-bf12f2

Archimedes (circa 287–212 BCE) has a right to a place among the great thinkers Raphael depicted in the fresco. From a visual point of view, it can be assumed (Joost-Gaugier, 2002) that Raphael projected Archimedes to sit in the place of the composition that Diogenes eventually occupies (just below Aristotle). Later, the artist removed this image, leaving only two mathematical protagonists, Pythagoras and Euclid.

In Ugo Panico da Carpi's chiaroscuro woodcut (Figure 7), Archimedes focuses on the rhombicuboctahedron, one of the thirteen Archimedean solids. Karpinski suggested that this is how Raphael depicted Archimedes in his lost preparatory drawing for the School of Athens. The artist chose this polyhedron as the most eloquent visual symbol of Archimedes' contribution to science (Karpinski, 2010).

Archimedean solids are semiregular convex polyhedra composed of two or more types of regular polygons meeting in identical vertices. They are so-called because Archimedes described them, even if we have only "second-hand" references to his writings on this topic. In Europe, the Archimedean polyhedra were rediscovered only in the Renaissance (Field, 1997).

Notably, the rhombicuboctahedron was first described by Luca Pacioli. In his book (Pacioli, 1509), Pacioli characterized this polyhedron as '*architecturally suggestive*', referring to a building tradition of the octagonal ground plan followed by pagan, Jewish, and Christian sacramental numerology of the number eight. Indeed, this polyhedron has eight-fold rotational symmetry and an octagonal prism in its central tier. Cross-sectioned, it corresponds to octagonal ground plans of many Renaissance religious buildings, for example, the Cathedral and Baptistery in Florence, most of Leonardo da Vinci's drawings of centrally planned churches (Xavier, 2008), and most importantly, the octagonal plan of St. Peter's Basilica, designed by Raphael's mentor Donato Bramante (1444 – 1514). The architecture of the building in the School of Athens was inspired by Bramante's design of St. Peter's Basilica. According to Vasari, Bramante helped Raphael with the architectural representation of the fresco (Janson, 2004).

Correspondingly, Caroline Karpinski hypothesized that Raphael chose the rhombicuboctahedron for Archimedes's iconography as a tribute to his mentor Donato Bramante (Karpinski, 2010).

On the other hand, the rhombicuboctahedron and the TEOB share the same properties: eight-fold rotational symmetry and an octagonal prism in its central tier. Accordingly, we can guess that Raphael chose the latter architecturally suggestive polyhedron for his drawing, shown in Fig. 7, for the same reason. This would be our answer to one of the questions formulated above.

## SUMMARY AND CONCLUSIONS

In the context of our previous research on geometrical forms of balance weights dated to the Early Islamic Period, we found one dirham weight with a complicated polyhedral shape in the Hecht Museum at the University of Haifa (Israel). We attributed it to a Truncated Elongated Octagonal Bipyramid (TEOB). Contrary to all Platonic solids and some Archimedean polyhedra, one cannot find any image or description of this polyhedron in the literature written by medieval Islamic mathematicians or early Renaissance European scientists. To our knowledge, the earliest Renaissance book containing the TEOB image is 'La Pratica di Prospettiva' published in 1596 by Florentine architect and perspective artist Lorenzo Sirigatti. We shortly described the book and discussed its importance as a practical manual for generations of artists and architects studying the perspective, in general, and light/shadow patterns of complex geometric forms, in particular. We mentioned Sirigatti's book's influence on Galileo Galilei's research and possible direct contact between them. We presented modern examples of TEOB-like lamps and lanterns Serendipitously, we revealed a drawing by Raphael with a lantern of a similar form.

# ACKNOWLEDGMENTS

I thank Dr. Vera Viana for bringing the image of the truncated elongated octagonal bipyramid in Lorenzo Sirigatti's "La pratica di prospettiva" to my attention and fruitful discussion about it. I also acknowledge assistance from Ms. Perry Harel, Registrar of the Hecht Museum at the University of Haifa, and virtual discussions with Prof. of Turin Polytechnic University Edoardo Piccoli and Dr. Lyubava Chistova from the State Hermitage Museum, St. Petersburg.